\newcommand{\ga}{\alpha}
\newcommand{\gb}{\beta}

\renewcommand{\gg}{\gamma}
\newcommand{\gd}{\delta}
\newcommand{\gn}{\nu}
\newcommand{\gm}{\mu}
\newcommand{\gl}{\lambda}
\newcommand{\gk}{\kappa}
\newcommand{\gr}{\rho}
\newcommand{\gx}{\xi}
\newcommand{\gf}{\phi}
\newcommand{\gq}{\theta}

\newcommand{\gt}{\tau}
\newcommand{\gs}{\sigma}
\newcommand{\go}{\omega}
\newcommand{\gp}{\pi}
\newcommand{\gh}{\eta}

\newcommand{\goy}{\omega_1}
\newcommand{\proof}{\noindent {\bf Proof.}\hspace{2mm}}

\def\rest{\mathord{\restriction}}
\newcommand{\force}{\Vdash} 
\newcommand{\til}[1]{\tilde{#1}}
\renewcommand{\aa}{\mbox{\sl aa\,}}
\newtheorem{theorem}{Theorem}
\newtheorem{corollary}[theorem]{Corollary}
\newtheorem{lemma}[theorem]{Lemma}
\newtheorem{proposition}[theorem]{Proposition}

\newtheorem{claim}{Claim}

\newtheorem{atheorem}{Theorem}
\newcommand{\se}{\subseteq}
\newcommand{\fin}{$\Box$\par\medskip}
\newcommand{\set}[2]{\{#1 \colon #2\}}
\newcommand{\ma}{{\cal A}}
\newcommand{\mb}{{\cal B}}
\newcommand{\mm}{{\cal M}}
\newcommand{\mn}{{\cal N}}
\newcommand{\mh}{{\cal H}}
\newcommand{\open}{\Bbb}

\newcommand{\oP}{{\open P}}
\newcommand{\oQ}{{\open Q}}

\newcommand{\oZ}{{\open Z}}
\newcommand{\EFG}[1]{{\cal G}_{#1}(\ma,\mb)}
\newcommand{\I}{\forall}
\newcommand{\II}{\exists}
\newcommand{\la}{\langle}
\newcommand{\ra}{\rangle}
\newcommand{\ol}[1]{\overline{#1}}

\newcommand{\cf}{\mbox{\rm cf}}
\newcommand{\dom}{\mbox{\rm dom}}

\documentstyle[12pt,amsfonts,amssymb]{article}
\author{by\\
\\
Alan Mekler\\
Department of Mathematics\\
Simon Fraser University\\
Burnaby, B.C.\\
Canada\\
\\
Saharon Shelah\\
Institute of Mathematics\\
Hebrew University\\
Jerusalem, Israel\\
\\
Jouko V\"a\"an\"anen\\
Department of Mathematics\\
University of Helsinki\\
Helsinki, Finland}
\title{The Ehrenfeucht-Fra\"\i ss\'e-game
of length~\(\goy\)}
\begin{document}
\maketitle

\begin{abstract}
Let \(\ma\) and \(\mb\) be two first order structures of the
same vocabulary. We shall consider the
{\em Ehrenfeucht-Fra\"\i ss\'e-game
of length \(\goy\) of \(\ma\) and \(\mb\)}
which we denote by \(\EFG{\goy}\). This game is like the
ordinary
Ehrenfeucht-Fra\"\i ss\'e-game of \(L_{\go\go}\) except that
there are \(\goy\) moves.
It is clear that \(\EFG{\goy}\) is determined if \(\ma\) and
\(\mb\)
are of cardinality \(\le\aleph_1\). We prove the following
results:
\begin{atheorem}
If V=L, then there are models \(\ma\) and \(\mb\) of
cardinality \(\aleph_2\) such that the game
\(\EFG{\goy}\) is non-determined.
\end{atheorem}
\begin{atheorem}
If it is consistent that there is a  measurable cardinal, then
it is consistent that \(\EFG{\goy}\) is determined for all
\(\ma\) and \(\mb\) of cardinality \(\le\aleph_2\).
\end{atheorem}
\begin{atheorem}
For any \(\gk\ge\aleph_3\)
there are \(\ma\) and \(\mb\) of cardinality \(\gk\)
such that the game \(\EFG{\goy}\) is
non-determined.
\end{atheorem}
\end{abstract}


\section{Introduction.}


Let \(\ma\) and \(\mb\) be two first order structures of the
same vocabulary \(L\). We denote the domains of \(\ma\)
and \(\mb\)
by \(A\) and \(B\) respectively.
All vocabularies are assumed to be relational.

The
{\em Ehrenfeucht-Fra\"\i ss\'e-game of length \(\gg\)
of \(\ma\) and \(\mb\)}
denoted  by \(\EFG{\gg}\) is defined as follows: There are
 two players
called \(\I\) and \(\II\). First \(\I\) plays \(x_0\) and
then \(\II\) plays
\(y_0\). After this
 \(\I\) plays \(x_1\), and \(\II\) plays \(y_1\), and
so on. If  \(\la(x_{\gb},y_{\gb}):\gb<\ga\ra\) has been
played
and \(\ga<\gg\), then  \(\I\) plays \(x_{\ga}\) after
which \(\II\)
plays \(y_{\ga}\). Eventually a sequence
\(\la(x_{\gb},y_{\gb}):\gb<\gg\ra\) has been played. The rules
of the game say that both players have to play elements
of \(A\cup B\). Moreover, if
\(\I\) plays his \(x_{\gb}\) in \(A\) (\(B\)), then \(\II\)
has to play
his \(y_{\gb}\) in \(B\) (\(A\)).
Thus the sequence \(\la(x_{\gb},y_{\gb}):\gb<\gg\ra\)
determines a relation \(\pi\se A\times B\).
Player \(\II\) wins this round of the game if \(\pi\) is
a partial isomorphism. Otherwise \(\I\) wins.
The notion of winning strategy is defined in the usual manner.
We say that a player {\em wins} \(\EFG{\gg}\) if he has
a winning
strategy in \(\EFG{\gg}\).

Recall that
\begin{eqnarray*}
\ma\equiv_{\go\go}\mb &\iff &\forall n<\go(\II\mbox{ wins
}\EFG{n})\\ \ma\equiv_{\infty\go}\mb &\iff &\II\mbox{ wins
}\EFG{\go}. \end{eqnarray*}
In particular, \(\EFG{\gg}\) is determined for \(\gg\le\go\).
The question, whether \(\EFG{\gg}\) is determined for
\(\gg>\go\),
is the subject of this paper. We shall concentrate on the
case  \(\gg=\goy\).

The notion
\begin{equation}
\label{IIwins}
\II\mbox{ wins }\EFG{\gg}
\end{equation}
can be viewed as a natural generalization of
\(\ma\equiv_{\infty\go}\mb\). The latter implies isomorphism
for countable models. Likewise (\ref{IIwins}) implies
isomorphism for models of cardinality \(|\gg|\):

\begin{proposition}
\label{first}
Suppose \(\ma\) and \(\mb\) have cardinality \(\le\gk\).
Then \(\EFG{\gk}\) is determined: \(\II\) wins if
\(\ma\cong\mb\), and
\(\I\) wins if \(\ma\not\cong\mb\).
\end{proposition}

\proof
If \(f:\ma\cong\mb\), then the winning strategy of \(\II\)
in \(\EFG{\gk}\)
is to play in such a way that the resulting \(\pi\) satisfies
\(\pi\se f\). On the other hand, if \(\ma\not\cong\mb\),
then the winning strategy of \(\I\) is to systematically
 enumerate
\(A\cup B\) so that the final \(\pi\) will satisfy
\(A= \dom (\pi)\) and \(B= rng(\pi)\). \fin

For models of arbitrary cardinality we have
the following simple but useful criterion of (\ref{IIwins}), namely in
the terminology of \cite{NS} that they are ``potentially isomorphic''.
We use
\(Col(\gl,\gk)\) to denote the notion of forcing which
collapses \(|\gl|\) to \(\gk\) (with conditions of cardinality less
than $\gk$).   

\begin{proposition}
\label{coll}
Suppose \(\ma\) and \(\mb\) have cardinality \(\le\gl\)
and \(\gk\) is regular.
Player \(\II\) wins \(\EFG{\gk}\) if and only if
\(\force_{Col(\gl,\gk)}\ma\cong\mb\).
\end{proposition}

\proof
Suppose \(\gt\) is a winning strategy of \(\II\)
in \(\EFG{\gk}\).
Since \(Col(\gl,\gk)\) is \(<\!\gk\)--closed,
\[\force_{Col(\gl,\gk)}``\gt\mbox{ is a winning
strategy of }\exists\mbox{ in }
G_{\gk}(\ma,\mb)".\]
Hence \(\force_{Col(\gl,\gk)}\ma\cong\mb\)  by
Proposition~\ref{first}.
Suppose then \(p\force \tilde{f}:\ma\cong\mb\)
for some \(p\in Col(\gl,\gk)\).
While the game \(\EFG{\gk}\) is played, \(\II\)
keeps extending the condition \(p\)
further and further. Suppose he has extended
\(p\) to \(q\) and \(\I\) has played \(x\in A\).
Then \(\II\) finds \(r\le q\) and \(y\in B\) with
\(r\force \tilde{f}(x)=y\). Using this simple strategy
\(\II\) wins. \fin


\begin{proposition}
\label{stable}
Suppose \(T\) is an $\go$-stable
 first order theory with NDOP.
Then \(\EFG{\goy}\) is determined for all
models \(\ma\) of \(T\) and all models \(\mb\).
\end{proposition}

\proof
Suppose \(\ma\) is a model of \(T\).
If \(\mb\) is not \(L_{\infty\omega_1}\)-equivalent to
\(\ma\), then \(\forall\) wins \(\EFG{\goy}\) easily.
So let us suppose \(\ma\equiv_{\infty\omega_1}\mb\).
We may assume \(A\) and \(B\) are of
cardinality \(\ge\aleph_1\). If we collapse
\(|A|\) and \(|B|\) to \(\aleph_1\), \(T\) will remain
$\go$-stable with NDOP, and \(\ma\) and \(\mb\)
will remain  \(L_{\infty\omega_1}\)-equivalent. So \(\ma\)
and \(\mb\) become isomorphic by
\cite[Chapter XIII, Section 1]{Shelah}.
Now Proposition 2 implies that \(\exists\) wins
\(\EFG{\goy}\). \fin

Hyttinen \cite{Hyt} showed that \(\EFG{\gg}\)
may be non-determined for all \(\gg\) with
\(\go<\gg<\goy\) and asked
whether
\(\EFG{\goy}\) may be non-determined. Our results show that
\(\EFG{\goy}\) may be non-determined for \(\ma\) and \(\mb\)
of cardinality \(\aleph_3\)
(Theorem~\ref{aleph_3}), but for models of cardinality
\(\aleph_2\) the answer is more complicated.

Let  \(F(\goy)\) be the free group of cardinality \(\aleph_1\).
Using the
combinatorial principle \(\Box_{\goy}\) we  construct
an abelian group \(G\)  of cardinality \(\aleph_2\)
such that \({\cal G}_{\goy}(F(\goy),G)\) is non-determined
(Theorem~\ref{aleph_2}).
 On the other hand,
we show that
starting with a model with  a measurable cardinal
one can build a forcing extension in which
\(\EFG{\goy}\) is determined for all models
\(\ma\) and \(\mb\) of cardinality \(\le\aleph_2\)
(Theorem~\ref{determined}).

Thus the free abelian group \(F(\goy)\) has the remarkable property
that the question
\begin{quote}
Is \({\cal G}_{\goy}(F(\goy),G)\) determined for all \(G\)?
\end{quote}
cannot be answered in ZFC alone. Proposition~\ref{stable}
shows that no model of an \(\aleph_1\)-categorical first order
theory  can have this property.

 We follow Jech \cite{Jech} in set theoretic notation.
We use \(S^m_n\) to denote the set
\(\{\ga<\go_m : \cf (\ga)=\go_n\}\).
Closed and unbounded sets are called cub sets.
A set  of ordinals is {\em \(\gl\)-closed}
if it is closed under supremums of ascending
\(\gl\)-sequences \(\la \ga_i : i<\gl\ra\) of
its elements. A subset of a cardinal is
{\em \(\gl\)-stationary} if it meets every
\(\gl\)-closed unbounded subset of the cardinal.
The closure of a set \(A\) of ordinals
in the order topology of ordinals
is denoted by \(\overline{A}\). The free
abelian group of cardinality \(\gk\)
is denoted by \(F(\gk)\).


\section{A non-determined
\({\cal G}_{\goy}(F(\goy),G)\) with \(G\) a group of cardinality
\(\aleph_2\).}


In this section we use \(\Box_{\goy}\) to construct
a group \(G\) of cardinality \(\aleph_2\)
such that the game \({\cal G}_{\goy}(F(\goy),G)\)
 is non-determined (Theorem~\ref{aleph_2}).
For background on almost free groups the reader is
referred to \cite{EM}. However, our presentation does not
depend
on special knowledge of almost free groups. All groups
below are assumed to be abelian.

By \(\Box_{\goy}\) we mean the  principle, which says that
there is a sequence
 \(\la C_{\ga} : \ga<\go_2, \ga=\cup\ga\ra\) such that
\begin{enumerate}
\item \(C_{\ga}\) is a cub subset of \(\ga\).
\item If \(\cf (\ga)=\go\), then \(|C_{\ga}|=\go\).
\item If \(\gg\) is a limit point of \(C_{\ga}\),
then \(C_{\gg}=C_{\ga}\cap\gg\).
\end{enumerate}
Recall that \(\Box_{\goy}\) follows from \(V=L\) by a result
of R. Jensen. For a sequence of sets \(C_{\ga}\)
as above we can let
\(E_{\gb}=\{\ga\in S^2_0 : \mbox{ the order type of }C_{\ga}
\mbox{ is }\gb\}\). For some \(\gb<\goy\) the set
\(E_{\gb}\) has to be stationary. Let us use \(E\)
to denote this \(E_{\gb}\). Then \(E\) is a so called
{\em non-reflecting} stationary set, i.e., if
\(\gg=\cup\gg\) then \(E\cap\gg\) is
non-stationary on \(\gg\). Indeed, then
some final segment \(D_{\gg}\) of
the set of limit points of \(C_{\gg}\)
is a cub subset of \(\gg\) disjoint from \(E\).
Moreover, \(\cf (\ga)=\go\) for all \(\ga\in E\).

\begin{theorem}
\label{aleph_2}
Assuming \(\Box_{\goy}\), there is a group \(G\)
 of cardinality \(\aleph_2\)
such that the game \({\cal G}_{\goy}(F(\goy),G)\) is non-determined.
\end{theorem}

\proof
Let \(\oZ^{\go_2}\) denote the direct product of \(\go_2\)
copies of the additive group \(\oZ\) of the integers.
Let \(x_{\ga}\) be the element of \(\oZ^{\go_2}\)
which is \(0\) on coordinates \(\not=\ga\)
and \(1\) on the coordinate \(\ga\). Let us fix
for each \(\gd\in S^2_0\) an ascending cofinal
sequence \(\gh_{\gd}:\go\to\gd\). For such \(\gd\), let
\[z_{\gd}=\sum_{n=0}^{\infty}2^nx_{\gh_{\gd}(n)}.\]
Let
\(\la C_{\ga} : \ga=\cup\ga<\go_2\ra\),
\(\la D_{\ga} : \ga=\cup\ga<\go_2\ra\) and
\(E=E_{\gb}\) be  obtained from \(\Box_{\goy}\)
as above. We are ready to define the groups
we need for the proof:
Let \(G\) be the smallest
pure subgroup of
\(\oZ^{\go_2}\) which contains \(x_{\ga}\) for \(\ga<\go_2\)
and \(z_{\gd}\) for \(\gd\in E\), let
\(G_{\ga}\) be the smallest
pure subgroup of
\(\oZ^{\go_2}\) which contains \(x_{\gg}\) for \(\gg<\ga\)
and \(z_{\gd}\) for \(\gd\in E\cap\ga\), let \(F\) \((=F(\go_2))\)
be the subgroup of
\(\oZ^{\go_2}\) generated freely by \(x_{\ga}\) for
\(\ga<\go_2\),
and finally, let \(F_{\ga}\)
be the subgroup of
\(\oZ^{\go_2}\) generated freely by \(x_{\gg}\) for
\(\gg<\ga\).


        The properties we shall want of $G_\ga$ are standard but for
the sake of completeness we shall sketch proofs. We need that each
$G_\ga$ is free and for any $\gb \notin E$ any free basis of $G_\gb$
can be extended to a free basis of $G_\ga$ for all $\ga > \gb$.

The proof is by induction on $\ga$. For limit ordinals we use the
fact that $E$ is non-reflecting. The case of successors of ordinals
not in $E$ is also easy. Assume now that $\gd \in E$ and the induction
hypothesis has been verified up to $\gd$. By the induction hypothesis
for any $\gb < \gd$ such that $\gb \notin E$, there is $n_0$ so that
$$G_\gd = G_\gb \oplus H \oplus K$$
where $K$ is the group freely
generated by $\set{x_{\eta_\gd(n)}}{n_0 \leq n}$ and $x_{\eta_\gd(m)}
\in G_\gb$ for all $m < n_0$. Then
        $$G_{\gd +1} = G_\gb \oplus H \oplus K'$$
where $K'$ is freely generated by
$$\set{\sum_{m=n}^{\infty}2^{m - n}x_{\gh_{\gd}(m)}}{n_0 \leq n}.$$

On the other hand, if \(\gd\in E\) and \(\{x_{\gh_{\gd}(n)} :
n<\go\}\se B\), where \(B\) is a subgroup of \(G\) such that $z_\gd
\notin B$,      
then \(G/B\) is non-free, as \(z_{\gd}+B\) is infinitely
divisible by \(2\) in \(G/B\).

\begin{claim}
\(\II\) does not  win
  \({\cal G}_{\goy}(F,G)\).
\end{claim}

 Suppose \(\gt\) is
a winning strategy of \(\II\). Let \(\ga\in E\)  such
that the pair \(( G_{\ga},F_{\ga})\) is closed
under the first \(\go\) moves of \(\gt\), that is,
if \(\I\) plays his first \(\go\) moves inside
\( G_{\ga}\cup F_{\ga}\), then \(\gt\) orders \(\II\)
to do the same. We shall  play \(G_{\goy}(F,A)\)
pointing out the moves of \(\I\) and letting
\(\gt\) determine the moves of \(\II\).
On his move number \(2n\)
\(\I\) plays the element \(x_{\gh_{\ga}(n)}\)
of \(G_{\ga}\). On his move number \(2n+1\)
\(\I\) plays some element of \(F_{\ga}\).
Player \(\I\) plays his moves in \(F_{\ga}\)
in such a way that during the first \(\go\) moves
eventually some countable direct summand \(K\) of
\(F_{\ga}\) as well as some countable \(B\se G_{\ga}\)
are enumerated. Let \(J\) be the smallest pure subgroup
of \(G\) containing \(B\cup\{z_{\ga}\}\). During the next
\(\go\) moves of \(G_{\goy}(F,A)\) player \(\I\)
enumerates \(J\) and \(\II\) responds by enumerating
some \(H\se F\). Since \(\gt\) is a winning strategy,
\(H\) has to be a subgroup of \(F\). But now \(H/K\)
is free, whereas \(J/B\) is non-free, so \(\I\) will win the game, a
contradiction.

\begin{claim} \(\I\) does not   win
  \({\cal G}_{\goy}(F,G)\).
\end{claim}

 Suppose \(\gt\) is
a winning strategy of \(\I\). If we were willing to use
CH, we could just take \(\ga\) of cofinality \(\goy\)
such that \((F_{\ga},G_{\ga})\) is closed under \(\gt\),
and derive a contradiction from the fact that
\(F_{\ga}\cong G_{\ga}\). However, since we do not want to
assume CH, we have to appeal to a longer argument.

Let \(\gk=(2^{\go})^{++}\). Let \(\mm\) be the expansion of
\(\la H(\gk),\in\ra\) obtained by adding the following
structure to it:
\begin{description}
\item[(H1)] The function \(\gd\mapsto\gh_{\gd}\).
\item[(H2)] The function \(\gd\mapsto z_{\gd}\).
\item[(H3)] The function \(\ga\mapsto C_{\ga}\).
\item[(H4)] A well-ordering \(<\) of the universe.
\item[(H5)] The winning strategy \(\gt\).
\end{description}
Let \(\mn=\la N,\in,\ldots \ra\) be an elementary submodel of
\(\mm\)
such that \(\goy\se N\) and \(N\cap\go_2\) is an ordinal
\(\ga\) of cofinality \(\goy\).

Let \(D_{\ga}=\{\gb_i : i<\goy\}\) in ascending order.
Since \(C_{\gb_i}=C_{\ga}\cap\gb_i\), every initial segment
of \(C_{\ga}\) is in \(N\). By elementaricity,
\(G_{\gb_i}\in N\) for all \(i<\goy\).
Let \(\gf\) be an isomorphism \(G_{\ga}\to F_{\ga}\)
obtained as follows: \(\gf\) restricted to \(G_{\gb_0}\)
is the \(<\)-least isomorphisms between the free
groups \(G_{\gb_0}\) and \(F_0\). If \(\gf\)
is defined on all \(G_{\gb_j}\), \(j<i\), then
\(\gf\) is defined on \(G_{\gb_i}\) as the \(<\)-least
extension of \(\bigcup_{j<i}\gf_{\gb_j}\) to an isomorphism
between \(G_{\gb_i}\) and \(F_i\). Recall that by our choice of
\(D_{\ga}\) \(G_{\gb_{i+1}}/G_{\gb_i}\) is free, so such
extensions really exist.

We derive a contradiction by showing that \(\II\) can play
\(\gf\) against \(\gt\) for the whole duration of the
game \(\EFG{\goy}\). To achieve this we have to show that,
when \(\II\) plays his canonical strategy based on \(\gf\)
the strategy \(\gt\) of \(\I\) directs \(\I\) to go on
playing elements which are in \(N\), that is, elements of
\(G_{\ga}\cup F_{\ga}\).

Suppose a sequence \(s=\la (x_{\gg},y_{\gg}) :
\gg<\gm\ra,\gm<\goy\), has been
played. It suffices to show that \(s\in N\). Choose \(\gb_i\)
so that the elements of \(s\) are in \(G_{\gb_i}\cup F_{\gb_i}\).
Now \(s\) is uniquely determined by \(\gf\rest G_{\gb_i}\)
and \(\gt\). Note that because
\(C_{\gb_i}=C_{\ga}\cap\gb_i\),
 \(\gf\rest G_{\gb_i}\) can be defined
inside \(N\) similarly as \(\gf\) was defined above,
using \(C_{\gb_i}\) instead of \(C_{\ga}\).
Thus \(s\in N\) and we are done.

We have proved that \({\cal G}_{\goy}(F,G)\) is nondetermined.
This clearly implies \({\cal G}_{\goy}(F(\goy),G)\)
is nondetermined.
 \fin

{\bf Remark.}
R. Jensen \cite[p. 286]{Jensen} showed that if \(\Box_{\goy}\)
fails, then \(\go_2\) is Mahlo in \(L\). Therefore, if
\(\EFG{\goy}\) is determined for all
almost free groups \(\ma\) and \(\mb\)
of cardinality \(\aleph_2\), then \(\go_2\) is Mahlo in \(L\).
If we start with \(\Box_{\gk}\), we get an almost free
group \(A\) of cardinality \(\gk^+\) such that
\({\cal G}_{\goy}(F(\goy),A)\) is nondetermined.


\section {\({\cal G}_{\goy}(F(\goy),G)\) can be determined for
all  \(G\).}


In this section all groups are assumed to be abelian.
It is easy to see that \(\II\) wins
\({\cal G}_{\goy}(F(\goy),G)\) for any uncountable
free group \(G\), so in this exposition \(F(\goy)\)
is a suitable representative of all free groups.
In the study of determinacy of
\({\cal G}_{\goy}(F(\goy),\ma)\)
for various \(\ma\) it suffices to study
groups \(\ma\), since for other \(\ma\)
 player \(\I\) easily wins the game.

Starting from a model with a Mahlo cardinal we
construct a forcing
extension in which \({\cal G}_{\goy}(F(\goy),G)\) is
determined, when
\(G\) is any
 group
 of cardinality \(\aleph_2\).
This can be extended to
groups \(G\) of any cardinality,
if we start with a supercompact
cardinal.

In the proof of the next results we shall make use of
{\em stationary logic} \(L(\aa)\). For the definition and basic
facts about \(L(\aa)\) the reader is referred to
\cite{BKM}. This logic has a new quantifier
 \(\aa s\) quantifying over variables \(s\)
ranging over countable subsets  of
the universe. A cub set of such \(s\) is any set which
 contains a
superset of any countable subset of the universe and
which is closed under unions of countable chains.
The semantics of aa \(s\) is defined as follows:
\[\aa s\gf(s,\ldots )\iff \gf(s,\ldots )
\mbox{ holds for a cub set of }s.\]
Note that a group of cardinality \(\aleph_1\) is free
if and only if it satisfies
\begin{equation}
\label{aa}
\aa s\:\aa s'(s\se s'\rightarrow s'/s
\mbox{ is free}).
\end{equation}

\begin{proposition}
\label{aa-char}
Let \(G\)  be a group. Then
the following conditions are equivalent:
\begin{description}
\item[(1)] \(\II\) wins
\({\cal G}_{\goy}(F(\goy),G)\).
\item[(2)]
\(G\) satisfies (\ref{aa}).
\item[(3)]   \(G\) is the union of
a continuous chain \(\la G_{\ga} : \ga<\go_2\ra\)
of free subgroups with \(G_{\ga+1}/G_{\ga}\) \(\aleph_1\)-free
for all \(\ga<\go_2\).
\end{description}
\end{proposition}

\proof (1) implies (2):
Suppose \(\II\) wins \({\cal G}_{\goy}(F(\goy),G)\).
By
 Proposition~\ref{coll} we have
\(\force_{Col(|G|,\goy)}``G\mbox{ is free}."\)
Using the countable completeness of
\(Col(|G|,\goy)\) it is now easy to construct
a cub set \(S\) of countable subgroups of \(G\)
such that if \(A\in S\) then for all
\(B\in S\) with \(A\se B\) we have \(B/A\) free.
Thus \(G\) satisfies (\ref{aa}).
(2) implies (3) quite trivially.
(3) implies (1):
Suppose  a continuous chain as in (3) exists.
If we collapse \(|G|\) to \(\aleph_1\), then
in the extension the chain has length \(<\go_2\).
Now we use Theorem 1 of \cite{Hill}:
\begin{quotation}
\noindent If a group \(A\) is the union of
a continuous chain of \(<\go_2\)
free subgroups \(\{A_{\ga}:\ga<\gg\}\)
of cardinality \(\le\aleph_1\) such that
 each \(A_{\ga+1}/A_{\ga}\) is \(\aleph_1\)-free,
then \(A\) is free.
\end{quotation}
Thus
\(G\) is free in the extension
and (1) follows from Proposition~\ref{coll}.
\fin


Let us consider the following principle:
\begin{description}
\item[\((*)\)]
For all stationary
\(E\se S^2_0\) and countable subsets
\(a_{\ga}\) of \(\ga\in E\) such that $a_\ga$ is cofinal in $\ga$ and
of order type $\go$
there is a closed \(C\se\go_2\)
of order type \(\goy\) such that
\(\{\ga\in E  : a_{\ga}\setminus C \mbox{ is finite}\}\) is
stationary in \(C\).
\end{description}

\begin{lemma}
The principle \((*)\) implies that
\({\cal G}_{\goy}(F(\goy),G)\)
is determined for all groups \(G\)
of cardinality
\(\aleph_2\).
\end{lemma}

\proof Suppose \(G\) is a group
of cardinality \(\aleph_2\).
We may assume the domain of \(G\) is \(\go_2\).
Let us assume \(G\) is \(\aleph_2\)-free,
as otherwise \(\I\) easily wins.
If we prove that \(G\) satisfies (\ref{aa}), then
Proposition~\ref{aa-char} implies that \(\II\)
wins \({\cal G}_{\goy}(F(\goy),G)\).

To prove (\ref{aa}), assume the contrary.
By Proposition~\ref{aa-char} we may assume
that  \(G\) can be expressed as the union of
a continuous chain \(\la G_{\ga} : \ga<\go_2\ra\)
of free groups with \(G_{\ga +1}/G_{\ga}\) non-\(\aleph_1\)-free
for  \(\ga\in E\), \(E\se\go_2\) stationary.
By Fodor's Lemma, we may assume \(E\se S^2_0\).  
Also we may assume
that for all $\ga$, every ordinal in $G_{\ga +1} \setminus G_\ga$ is
greater than
every ordinal in $G_\ga$. Finally by intersecting with a closed
unbounded set we may assume that for all $\ga \in E$ the set
underlying $G_\ga$ is $\ga$.
Choose for each
\(\ga\in E\) some countable subgroup
\(b_{\ga}\) of \(G_{\ga +1}\) with
\(b_{\ga}+G_{\ga}/G_{\ga}\) non-free.
Let \(c_{\ga}=b_{\ga}\cap G_{\ga}\). 
We will choose $a_\ga$ so that any final segment generates a subgroup
containing $c_\ga$. Enumerate $c_\ga$ as $\set{g_n}{n < \go}$ such
that each element is enumerated infintely often. Choose an
increasing sequence $(\ga_n\colon n < \go)$ cofinal in $\ga$ so that
for all $n$, $g_n \in G_{\ga_n}$. Finally, for each $n$, choose $h_n
\in G_{\ga_n +1} \setminus G_{\ga_n}$. Let $a_\ga = \set{h_n}{n <
\go} \cup \set{h_n + g_n}{n < \go}$. It is now easy to check that
$a_\ga$ is a sequence of order type $\go$ which is cofinal in $\ga$ and any
subgroup of $G$ which contains all but finitely many of the elements
of $a_\ga$ contains $c_\ga$.

        By \((*)\) there is a continuous
\(C\) of order type \(\goy\) such that
\(\{\ga\in C : a_{\ga}\setminus C \mbox{ is finite}\}\) is stationary
in \(C\).
Let \(D=\la C \bigcup \sum_{\ga\in C}b_{\ga}\ra\).
Since \(|D|\le\aleph_1\), \(D\) is free.

        For any $\ga\in C$, let
$$D_\ga = \la (C \cap \ga) \bigcup (\sum_{\gb \in (C \cap \ga)} b_\gb)\ra.$$
Note that $D = \bigcup_{\ga \in C} D_\ga$, each $D_\ga$ is countable
and for limit point $\gd$ of $C$, $D_\gd = \bigcup_{\ga \in (C \cap
\gd)} D_\ga$. Hence
there is an \(\ga\in C\cap E\) such that
\(a_{\ga}\setminus C\) is finite
and \(D/D_{\ga}\) is free.
Hence \(b_{\ga} + D_\ga/D_\ga\) is free. But
\[b_{\ga} + D_\ga/D_\ga
\cong b_\ga/b_{\ga}\cap D_{\ga}=b_{\ga}/b_\ga \cap G_{\ga},\]
which is not free, a contradiction.
\fin

For the next theorem we need a
lemma from \cite{GiSh}. A proof is
included for the convenience of the reader.

\begin{lemma}
\label{stat}
\cite{GiSh} Suppose $\lambda$ is a regular cardinal
and ${\open Q}$ is a
notion of
forcing which satisfies the $\lambda$-c.c.
Suppose $\cal I$ is a
normal $\gl$-complete ideal on $\gl$
and \({\cal I}^+=\{S\se\gl : S\not\in {\cal I}\}\).
For all
sets $S\in{\cal I}^+$ and sequences of
conditions $\la p_\alpha \colon
\alpha \in S\ra$, there is a set
$C$ with \(\gl\setminus C\in{\cal I}\) so that for
all $\alpha \in C \cap S$,
\[p_\alpha \force_{\oQ}
``\set{\gb}{p_\gb \in \til{G}}\in {\cal J}^+ \mbox{, where }
{\cal J}\mbox{ is
the ideal generated by $\cal I$}".\]
\end{lemma}

\proof Suppose the lemma is false. So there is an
${\cal I}$-positive
set $S' \se S$ such that for all $\ga \in S'$ there
is an extension
$r_\ga$ of $p_\ga$ and a set $I_\ga \in \cal I$
(note: $I_\ga$ is in
the ground model) so that
\[r_\ga \force \{\gb : p_{\gb}\in\tilde{G}\}
\se I_\ga.\]
Let $I$ be the diagonal union of $\set{I_\ga}{\ga \in S'}$.

        Suppose now that $\ga < \gb$ and $\ga, \gb \in
(S' \setminus
I)$. Since $\gb \notin I$, $r_\ga \force
p_{\gb}\not\in\tilde{G}$.
Hence $r_\ga \force r_\gb \not\in
\tilde{G}$. So $r_\ga, r_\gb$ are incompatible. Hence
$\set{r_\ga}{\ga
\in S' \setminus I}$ is an antichain which, since $S'$ is $\cal
I$-positive, is of cardinality $\gl$. This is a contradiction.
\fin

\begin{theorem}
\label{Mahlo}
Assuming the consistency of a Mahlo cardinal, it is
consistent that \((*)\) holds and hence that
\({\cal G}_{\goy}(F(\goy),G)\)
is determined for all  groups
\(G\) of  cardinality
\(\aleph_2\).
\end{theorem}

\proof
By a result of Harrington and Shelah \cite{HaSh} we may start
with a Mahlo cardinal \(\gk\) in which every stationary set
of cofinality \(\go\) reflects,
that is, if \(S\se\gk\) is stationary,
and \(\cf (\ga)=\go\) for \(\ga\in S\), then \(S\cap\gl\)
is stationary in \(\gl\)
for some inaccessible \(\gl<\gk\).

For any inaccessible \(\gl\)
let \(\oP_{\gl}\) be the Levy-forcing
for collapsing \(\gl\) to \(\go_2\).
The conditions of \(\oP_{\gl}\) are countable functions
\(f:\gl\times\goy\to\gl\) such that \(f(\ga,\gb)<\ga\) for all
\(\ga\) and \(\gb\) and each $f$ is increasing and continuous in the
second coordinate. 
It is well-known that \(\oP_{\gl}\)
is countably closed and satisfies the \(\gl\)-chain condition
\cite[p. 191]{Jech}.

Let \(\oP=\oP_{\gk}\).
Suppose \(p\in\oP\)
and \[p\force``\til{E}\se S^2_0 \mbox{ is stationary and }
\forall\ga\in\til{E}(\til{a}_{\ga}\se\ga
\mbox{ is cofinal in $\ga$ and of order type $\go$)"}.\]
Let \[S=\{\ga<\gk : \exists q\le p(q\force \ga\in\til{E})\}.\]
For any \(\ga\in S\) let \(p_{\ga}\le p\)
such that \(p_{\ga}\force\ga\in\til{E}\). Since \(\oP\) is countably
closed, we
can additionally require that for some  countable
\(a_{\ga}\se\ga\) we have
\(p_{\ga}\force\til{a}_{\ga}=a_{\ga}\).

The set \(S\) is stationary in
\(\gk\), for if \(C\se\gk\) is cub, then
\(p\force C\cap\til{E}\not=\emptyset\), whence
\(C\cap S\not=\emptyset\). Also
\(\cf (\ga)=\go\) for \(\ga\in S\).
Let \(\gl\) be inaccessible such that \(S\cap\gl\)
is stationary in \(\gl\). We may choose \(\gl\)
in such a way that \(\ga\in S\cap\gl\) implies
\(p_{\ga}\in\oP_{\gl}\).
By Lemma~\ref{stat} there is a \(\gd\in S\cap\gl\)
such that
\[p_{\gd}\force_{\oP_{\gl}}``\til{E_1}=\{\ga<\gl :
p_{\ga}\in\til{G}\}
\mbox{ is stationary.}"\]

Let  \(\oQ\) be the set of
conditions \(f\in \oP\) with
\(\dom (f)\se(\gk\setminus\gl)\times\goy\). Note that
\(\oP\cong\oP_{\gl}\otimes\oQ\).
Let \(G\) be \(\oP\)-generic containing \(p_{\gd}\) and
\(G_{\gl}=G\cap\oP_{\gl}\)
for any inaccessible \(\gl\le\gk\). Then \(G_{\gl}\) is
\(\oP_{\gl}\)-generic and  \(\go_2\) of  \(V[G_{\gl}]\)
is \(\gl\).
Let us work now in \(V[G_{\gl}]\). Thus \(\gl\) is the current
\(\go_2\), \(E_1=\{\ga<\gl : p_{\ga}\in G_{\gl}\}\) is
stationary,
and we have the countable sets \(a_{\ga}\se\ga\) for
\(\ga\in E_1\).
Since \(\oQ\) collapses \(\gl\) there is a name \(\til{f}\)
such that
\[\force_{\oQ}``\til{f}:\goy\to\gl\mbox{ is continuous
and cofinal}."\]
More precisely $\til{f}$ is the name for the function $f$ defined by
$f(\ga) = \gb$ if and only if there is some $g \in G$ so that $g(\gl,
\ga) = \gb$.  
Let \(\til{C}\) denote the range of \(\til{f}\).
We shall prove the following statement:
\medskip

\noindent{\bf Claim:}
\(\force_{\oQ}
\{\ga<\gl  : a_{\ga}\setminus \til{C} \mbox{ is finite }\}\mbox{ is
stationary in }
\til{C}.
\)
\medskip


Suppose \(q\in\oQ\) so that
\(q\force``\til{D}\se\goy\mbox{ is a cub}."\)
Let \(\mm\) be an appropriate expansion of
\(\la H(\gk),\in\ra\) and \(\la \mn_i : i<\gl\ra\),
\(\mn_i=\la N_i,\in,\ldots \ra\),
a sequence of elementary submodels of \(\mm\)
such that:
\begin{description}
\item[(i)] Everything relevant is in \(N_0\).
\item[(ii)] If \(\ga_i=N_i\cap\gl\), then
\(\ga_i<\ga_j\) for \(i<j<\gl\).
\item[(iii)] \(N_{i+1}\) is closed under countable sequences.
\item[(iv)] \(|N_i|=\goy\).
\item[(v)] \(N_i=\bigcup_{j<i}N_j\) for \(i\) a limit ordinal.
\end{description}

        Choose \(\gg=\ga_i\in E_1\) and let \(\la i_n :n<\go\ra\) be
a sequence of successor ordinals such that
\(\gg=\sup \{\ga_{i_n}:n<\go\}\).
Let \(q_0\le q\) and \(\gb_0 \in\goy\) such that
\(q_0 , \gb_0 \in N_{i_0}\),
\[q_0\force ``\gb_0\in\til{D}\mbox{''}\]
 and $q_0$ decides the value of $\til{f}''\gb_0$ (which will by
elementaricity necessarilly be a subset of $\ga_{i_0}$).

If \(q_n\) and \(\gb_n\) are defined we choose
\(q_{n+1}\le q_n\)
and \(\gb_{n+1}\in\goy\) such that
\(q_{n+1}, \gb_{n+1}\in N_{i_{n+1}}\),
\[q_{n+1}\force ``\gb_{n+1}\in\til{D}\mbox{ and }
a_{\gg}\cap(\ga_{i_{n+1}}\setminus \ga_{i_n})\se\til{f}\mbox{''}\gb_{n+1}\se
\ga_{i_{n+1}}"\]
and
\(q_{n+1}\) decides \(\til{f}\mbox{''}\gb_{n+1}\).
Finally, let \(q_{\go}=\bigcup\{q_n:n<\go\}\)
and \(\gb=\bigcup\{\gb_n:n<\go\}\). Then
\[q_{\go}\force``\gb\in\til{D}
\mbox{ and }
a_{\gg}\setminus\til{f}\mbox{''}\gb\mbox{ is finite}."\]
The claim, and thereby the theorem, is proved. \fin

\begin{corollary}
The statement that
\(\EFG{\goy}\)
is determined for every structure \(\ma\) of cardinality
\(\aleph_2\)
and every uncountable free  group \(\mb\), is
equiconsistent with the existence of a Mahlo cardinal.
\end{corollary}

\noindent{\bf Remark.} If \({\cal G}_{\goy}(A,F(\go_1))\) is
determined for all groups
\(A\)
of  cardinality \(\gk^+\), \(\gk\) singular,
then \(\Box_{\gk}\) fails. This implies
that the Covering Lemma fails for the Core Model,
whence there is an inner model for a measurable
cardinal. This shows that the conclusion of Theorem~\ref{Mahlo}
cannot be strengthened to arbitrary \(G\).
However, by starting with a larger cardinal we can
make this extension:

\begin{theorem}
\label{succ}
Assuming the consistency of a supercompact cardinal, it is
consistent that
\({\cal G}_{\goy}(F(\goy),G)\)
is determined for all  groups
\(G\).
\end{theorem}

\proof
Let us assume that the stationary logic \(L_{\go_1\go}(aa)\)
has the  L\"owenheim-Skolem property down to \(\aleph_1\).
This assumption is
consistent relative to the consistency of a supercompact
cardinal \cite{BD}. Let \(G\) be an arbitrary \(\aleph_2\)-free
group. Let \(H\) be an \(L(aa)\)-elementary
submodel of \(G\) of cardinality \(\aleph_1\).
Thus \(H\) is a free group.
The group \(H\) satisfies the sentence (\ref{aa}),
whence so does \(G\). Now the claim follows from
Proposition~\ref{aa-char}. \fin

\begin{corollary}
Assuming the consistency of a supercompact cardinal, it is
consistent that
\(\EFG{\goy}\)
is determined for every structure \(\ma\)
and every uncountable free  group \(\mb\).
\end{corollary}


\section {\(\EFG{\goy}\)
can be determined for all
\(\ma\) and \(\mb\) of  cardinality
\(\aleph_2\).}

We prove the  consistency of the statement
that \(\EFG{\goy}\) is determined for all
\(\ma\) and \(\mb\) of  cardinality
\(\le\aleph_2\) assuming the consistency
of a measurable cardinal. Actually we make
use of an assumption that we call \(I^*(\go)\) concerning
stationary subsets of \(\go_2\). This assumption
is known to imply that \(\go_2\) is measurable in
an inner model. It follows from the previous section
that some large cardinal axioms are needed
to prove the stated determinacy.

Let \(I^*(\go)\) be the following assumption about
\(\goy\)-stationary subsets of \(\go_2\):

\begin{description}
\item[\(I^*(\go)\)] Let \({\cal I}\) be
the  \(\goy\)-nonstationary ideal \(NS_{\goy}\) on
 \(\go_2\). Then \({\cal I}^+\) has a \(\gs\)-closed dense subset \(K\).
\end{description}

\noindent
Hodges and Shelah \cite{HS} define a  principle \(I(\go)\),
which is like
\(I^*(\go)\) except that \({\cal I}\) is not assumed to be the
\(\goy\)-nonstationary ideal. They use \(I(\go)\) to prove the
determinacy of an Ehrenfeucht-Fra\"\i ss\'e-game played on
several boards simulataneously.

Note that \(I^*(\go)\) implies \({\cal I}\) is precipitous, so the
consistency
of \(I^*(\go)\) implies the consistency of a measurable cardinal
\cite{JMMP}.

\begin{theorem}
\label{I}
(\cite{JMMP})The assumption \(I^*(\go)\) is consistent
 relative to
the consistency of a measurable cardinal.
\end{theorem}

We shall consider models
\(\ma,\mb\) of cardinality \(\aleph_2\), so we may as well
assume they have \(\go_2\) as universe. For such \(\ma\)
and \(\ga<\go_2\) we let \(\ma_{\ga}\) denote the structure
\(\ma\cap\ga\).
Similarly \(\mb_{\ga}\).

\begin{lemma}
\label{stationary}
Suppose \(\ma\) and \(\mb\) are structures of
cardinality \(\aleph_2\).
If \(\I\) does not have a winning strategy in
\(\EFG{\goy}\), then
\[S=\{\ga : \ma_{\ga}\cong\mb_{\ga}\}\]
is \(\goy\)-stationary.
\end{lemma}

\proof
Let \(C\se\go_2\) be \(\goy\)-closed and unbounded. Suppose
\(S\cap C=\emptyset\).
We derive a contradiction by describing a winning strategy of
\(\I\):
Let \(\gp:\go_1\to\go_1\times\goy\times 2\) be onto with
\(\ga,\gb,d\le\gp(\ga,\gb,d)\) for all \(\ga,\gb<\goy\) and
\(d<2\). If \(\ga<\go_2\), let \(\gq_{\ga}:\goy\to\ga\) be onto.
Suppose the sequence \(\la (x_i,y_i) : i<\ga\ra\) has been
played.
Here \(x_i\) denotes a move of \(\I\) and \(y_i\) a move
of \(\II\). During the game \(\I\) has built
an ascending sequence \(\{c_i: i<\ga\}\) of elements of \(C\).
Now he lets \(c_{\ga}\) be the smallest element of \(C\)
greater
than all the elements \(x_i,y_i,i<\ga\). Suppose
\(\gp(\ga)=(i,\gg,d)\). Now \(\I\) will play
\(\gq_{c_i}(\gg)\) as an element of \(\ma\), if \(d=0\), and
as an element of \(\mb\) if \(d=1\).

After all \(\goy\) moves of \(\EFG{\goy}\) have been played,
some \(\ma_{\ga}\) and \(\mb_{\ga}\), where \(\ga\in C\),
have been
enumerated. Since \(\ga\not\in S\), \(\I\) has won the game.
\fin

\begin{theorem}
\label{determined}
Assume \(I^*(\go)\). The game \(\EFG{\goy}\) is determined for
all \(\ma\) and \(\mb\) of cardinality \(\le\aleph_2\).
\end{theorem}

\proof
Suppose \(\I\) does not have a winning strategy.
By  Lemma~\ref{stationary} the set \(S=\{\ga :
\ma_{\ga}\cong\mb_{\ga}\}\)
is \(\goy\)-stationary. Let \(I\) and \(K\) be as in
\(I^*(\go)\).
If \(\ga\in S\), let \(h_{\ga}:\ma_{\ga}\cong\mb_{\ga}\).
We describe
a winning strategy of \(\II\).
The idea of this strategy is that \(\II\) lets
the isomorphisms \(h_{\ga}\) determine his moves.
Of course, different \(h_{\ga}\) may give different
information to \(\II\), so he has to decide which
\(h_{\ga}\) to follow. The key point is that \(\II\)
lets some \(h_{\ga}\) determine his move only if there
are stationarily many other \(h_{\gb}\) that agree
with \(h_{\ga}\) on this move.

Suppose the sequence \(\la (x_i,y_i) : i<\ga\ra\) has been
played.
Again \(x_i\) denotes a move of \(\I\) and \(y_i\) a move
of \(\II\).
Suppose \(\I\) plays next
\(x_{\ga}\) and this is (say) in \(A\).
During the game \(\II\) has built
a  descending sequence \(\{S_i: i<\ga\}\) of elements of \(K\)
with \(S_0\se S\).
The point of the sets \(S_i\) is
that \(\II\) has taken care that
for all \(i<\ga\) and \(\gb\in S_i\)
we have
\(y_i=h_{\gb}(x_i)\)
or
\(x_i=h_{\gb}(y_i)\)
depending on whether \(\I\) played \(x_i\) in
\(A\) or \(B\).
Now \(\II\) lets
\(S'_{\ga}\se\bigcap_{i<\ga}S_i\) so that  \(S'_{\ga}\in K\)
and
\(\forall i\in S'_{\ga}(x_{\ga}<i)\).
 For each \(i\in S'_{\ga}\)
we have \(h_i(x_{\ga})<i\). By normality, there are an
\(S_{\ga}\se S'_{\ga}\) in \(K\) and a \(y_{\ga}\) such that
\(h_i(x_{\ga})=y_{\ga}\) for all \(i\in S_{\ga}\). This element
\(y_{\ga}\) is the next move of \(\II\). Using this
strategy \(\II\)
wins. \fin


\section{A non-determined
\(\EFG{\goy}\) with \(\ma\) and \(\mb\) of cardinality
\(\aleph_3\).}


We construct directly in ZFC two models
\(\ma\) and \(\mb\) of cardinality
\(\aleph_3\) with
\(\EFG{\goy}\)  non-determined.
It readily follows that
such models exist in all cardinalities
\(\ge\aleph_3\). The construction uses a
square-like principle
(Lemma~\ref{comb}), which is provable
in ZFC.

\begin{lemma}
\label{BM}\cite{Sh89,ShBook}
There is a stationary \(X\se S^3_1\) and a sequence
\(\la D_{\ga} : \ga\in X\ra\) such that
\begin{enumerate}
\item  \(D_{\ga}\) is a cub subset of \(\ga\) for all
\(\ga\in X\).
\item  The order type of \(D_{\ga}\) is \(\goy\).
\item If \(\ga,\gb\in X\) and \(\gg<min\{\ga,\gb\}\) is
a limit of both \(D_{\ga}\) and \(D_{\gb}\),
then \(D_{\ga}\cap\gg = D_{\gb}\cap\gg\).
\item If \(\gg\in D_{\ga}\), then \(\gg\) is a limit point
of \(D_{\ga}\) if and only if \(\gg\) is a limit ordinal.
\end{enumerate}
\end{lemma}

\proof
We shall sketch,
for completeness,
 a proof of this given by Burke and
Magidor~\cite[Lemma
7.7]{BM}.

Let \(<^*\) be a well-ordering of \(H(\go_3)\).
For each \(\ga\in S^3_1\), let
\(\la N^{\ga}_{\gd} : \gd<\go_2\ra\)
be a continuously increasing chain of elementary submodels
of \(\la H(\go_3),\in,<^*\ra\) such that
\begin{description}
\item[(N1)] \((\goy+1)\cup\{\go_2,\ga\}\se N^{\ga}_0\).
\item[(N2)] \(|N^{\ga}_{\gd}|\le\goy\).
\item[(N3)] \(N^{\ga}_{\gd}\cap\go_2\in\go_2\).
\item[(N4)] \(\overline{N^{\ga}_{\gd}\cap\go_3}\in
N^{\ga}_{\gd+1}\).
\end{description}
Let \(A^{\ga}_{\gd}=N^{\ga}_{\gd}\cap\ga\) for each
\(\ga\in S^3_1\).
Since, \(\ga\in N^{\ga}_{\gd}\), \(A^{\ga}_{\gd}\) is cofinal
in \(\ga\).
Let \(X\se S^3_1\) be stationary such that for some
\(\gd,\gr<\go_2\)
and for all \(\ga\in X\) we have
\begin{enumerate}
\item \(\gd\)= least ordinal of cofinality \(\goy\) with
\(N^{\ga}_{\gd}\cap\go_2=\gd\).
\item The order type of \(\overline{A^{\ga}_{\gd}}\) is
\(\gr+1\).
\end{enumerate}
Let \(f:\goy\to\gr\)
be cofinal and continuous.
Let \(g:\gr+1\cong\overline{A^{\ga}_{\gd}}\)
such that \(gf\) maps successors to successors.
Let \(D_{\ga}\) be the image of \(\goy\) under \(gf\).
\fin

\begin{lemma}
\label{comb}
There are sets \(S,T\) and \(C_{\ga}\) for \(\ga\in S\) such
that
the following hold:
\begin{enumerate}
\item \(S\se S_0^3\cup S_1^3\) and \(S\cap S^3_1\) is stationary.
\item \(T\se S^3_0\) is stationary and \(S\cap T=\emptyset\).
\item If \(\ga\in S\), then
\(C_{\ga}\se\ga\cap S\) is closed
 and of order-type \(\le\goy\).
\item If \(\ga\in S\) and
\(\gb\in C_{\ga}\), then \(C_{\gb}=C_{\ga}\cap\gb\).
\item If \(\ga\in S\cap S^3_1\), then \(C_{\ga}\) is cub on
\(\ga\).
\end{enumerate}
\end{lemma}

\proof
 Let \(S\) and \(\la D_{\ga} : \ga\in S\ra\) be as in
Lemma~\ref{BM}.
Let \(S'=X\cup Y\), where \(Y\) consists of ordinals which are
limit points \(<\ga\) of some \(D_{\ga},\ga\in X\).
If \(\ga\in X\),
we let \(C_{\ga}\) be the set of limit points \(<\ga\)
of \(D_{\ga}\).
If \(\ga\in Y\), we let \(C_{\ga}\) be the set of limit points
\(<\ga\) of \(D_{\gb}\cap\ga\), where \(\gb>\ga\) is chosen
arbitrarily from \(X\).

Now claims 1,3,4 and 6 are clearly satisfied.

Let \(S^3_0=\bigcup_{i<\go_2}T_i\) where the \(T_i\) are
disjoint
stationary sets.
Since \(|\overline{C_{\ga}}|\le\goy\), there is
\(i_{\ga}<\go_2\) such that \(i\ge i_{\ga}\) implies
\(\overline{C_i}\cap T_i=\emptyset.\) Let \(S''\se S'\) be
stationary such
that \(\ga\in S''\) implies \(i_{\ga}\) is constant \(i\).
Let \(T=T_i\). Finally, let
\(S=S''\cup\bigcup\{C_{\ga} : \ga\in S''\}\).
Claim 2 is satisfied, and the Lemma is proved.
\fin

\begin{theorem}
\label{aleph_3}
There are structures \(\ma\) and \(\mb\) of cardinality
\(\aleph_3\) with one binary predicate such that the game
\(\EFG{\goy}\) is non-determined.
\end{theorem}

\proof
Let \(S,T\) and \(\la C_{\ga} : \ga\in S\ra\)
be as in Lemma~\ref{comb}.
We shall
construct a sequence \(\{M_{\ga} : \ga<\go_3\}\)
of sets and a sequence \(\{G_{\ga} : \ga\in S\}\) of
functions such that the conditions
(M1)--(M6) below hold.
Let \(W_{\ga}\) be the set of all mappings
\[G^{d_0}_{\gg_0}\ldots G^{d_n}_{\gg_n},\]
where
\(\gg_0,\ldots ,\gg_n\in S\cap\ga\), \(d_i\in\{-1,1\}\),
\(G^{1}_{\gg}\) means \(G_{\gg}\)
and \(G^{-1}_{\gg}\) means
the inverse of \(G_{\gg}\).
Let \(W=W_{\go_3}\). (Note that $W$ consists of a set of partial
functions.)

The conditions on the \(M_{\ga}\)'s and the \(G_{\ga}\)'s
are:
\begin{description}
\item[(M1)] \(M_{\ga}\subseteq M_{\gb}\) if \(\ga<\gb\),
and \(M_{\ga}\subset M_{\ga+1}\) if  \(\ga\in S\).
\item[(M2)] \(M_{\gn}=\bigcup_{\ga<\gn}M_{\ga}\) for limit
\(\gn\).
\item[(M3)] \(G_{\ga}\) is a bijection of \(M_{\ga+1}\) for
\(\ga\in S\).
\item[(M4)] If \(\gb\in S\) and \(\ga\in C_{\gb}\), then
\(G_{\ga}\se G_{\gb}\).
\item[(M5)] If for some $\gb$, $G_\gb(a) = b$ and for some $w \in W$,
$w(a) = b$, then there is some $\gamma$ so that $w \se G_\gamma$.
Furthermore if $\gb$ is the minimum ordinal so that $G_\gb(a) = b$
then $\gamma = \gb$ or $\gb \in C_\gamma$.
\end{description}

In order to construct the set \(M=\bigcup_{\ga<\go_3}M_{\ga}\) and the
mappings  \(G_{\ga}\) we define an oriented graph with \(M\) as the
set of vertices.  We use the terminology of Serre~\cite{S} for
graph-theoretic notions.  If \(x\) is an edge, the origin of \(x\) is
denoted by \(o(x)\) and the terminus by \(t(x)\).  Our graph has an
inverse edge \(\ol{x}\) for each edge \(x\). Thus \(o(\ol{x})=t(x)\)
and \(t(\ol{x})=o(x)\).  Some edges are called {\em positive}, the
rest are called {\em negative}. An edge is positive if and only if its
inverse is negative. For each edge \(x\) of \(M\) there is a set
\(L(x)\) of labels. The set of possible labels for positive edges is
\(\{g_{\ga} : \ga<\go_3\}.\) The negative edges can have elements of
\(\{g^{-1}_{\ga} : \ga<\go_3\}\) as labels. The labels are assumed to
be given in such a way that a positive edge gets \(g_{\ga}\) as a
label if and only if its inverse gets the label \(g^{-1}_{\ga}\).
During the construction the sets of labels will be extended step by
step.


        The construction is analogous to building an acyclic graph on
which a group acts freely. The graph then turns out to be the Cayley
graph of the group. The labelled graph we will build will be the
``Cayley graph'' of $W$ which will be as free as possible given
(M1)--(M4). Condition (M5) is a consequence of the freeness of the
construction.

Let us suppose the sets \(M_{\gb},\gb<\ga\), of vertices have
been
defined. Let \(M_{<\ga}=\bigcup_{\gb<\ga}M_{\ga}\).
Some vertices in \(M_{<\ga}\) have edges between them
and a set \(L(x)\) of labels has been assigned to each such
edge \(x\).

If \(\ga\) is a limit ordinal, we let \(M_{\ga}=M_{<\ga}\).
So let us assume \(\ga=\gb+1\).
If \(\gb\not\in S\), \(M_{\ga}=M_{\gb}\). So let
us assume \(\gb\in S\). Let \(\gg=\sup(C_{\gb})\). 
Notice that since $S$ consists entirely of limit ordinals and $C_\gb
\se S$, either $\gg = \gb$ or $\gg +1 < \gb$.

\medskip
\noindent{\bf Case 1.} \(\gg=\gb\): We extend
\(M_{\gb}\)
to \(M_{\ga}\) by adding new vertices \(\{P_z : z\in \oZ\}\)
and for each \(z\in\oZ\) a positive edge
\(x^{P_z}_{\ga}\) with \(o(x^{P_z}_{\ga})=P_z\)
and \(t(x^{P_z}_{\ga})=P_{z+1}\).
We also let \(L(x^{P_z}_{\ga})= \{g_{\gb}\} \cup \set{g_\gd}{\gb \in C_\gd}\).
\medskip

\noindent{\bf Case 2.} \(\gg+1<\gb\): We extend \(M_{\gb}\)
to \(M_{\ga}\) by adding new vertices
\(\{P'_z : z\in \oZ \setminus \{0\}\}\) for each
\(P\in M_{\gb}\setminus M_{\gg+1}\). For notational convenience let
$P'_0 = P$. Now we add for each
\(P\in M_{\gb}\setminus M_{\gg+1}\)
new edges as follows.
For each \(z\in\oZ\) we add a positive
edge \(x^{P'_z}_{\ga}\) with
\[
o(x^{P'_z}_{\ga})       = P'_z,
t(x^{P'_z}_{\ga}) = P'_{z+1},
L(x^{P'_z}_{\ga}) = \{g_{\gb}\} \cup \set{g_\gd}{\gb \in C_\gd}
\]
This determines completely the inverse of \(x^{P'_z}_{\ga}\).

\medskip
This ends the construction of the graph.
In the construction each vertex \(P\) in
\(M_{\ga+1}\),\(\ga\in S\), is made
the origin of a unique edge \(x^P_{\ga}\) with
\(g_{\ga}\in L(x^P_{\ga})\). We define
\[G_{\ga}(P)=t(x^P_{\ga}).\]

The construction of the sets \(M_{\ga}\) and
the mappings \(G_{\ga}\) is now completed.
It follows immediately from the construction that
each \(G_{\ga}\),
\(\ga\in S\), is a bijection of \(M_{\ga+1}\).
So (M1)--(M3) hold. (M4) holds, because
\(g_{\ga}\) is added to the labels of any
edge with \(g_{\gb}\), where \(\gb\in C_{\ga}\), as a label.
Finally, (M5) is a consequence of the fact that the
graph is circuit-free.

Let us fix \(a_0\in M_1\)
and \(b_0=G_{\gb_0}(a_0)\),
where \(\gb_0\in C_{\ga}\)
for all \(\ga\in S\). Note that we may assume, without loss of
generality, the existence of such a \(\gb_0\).

If \(a_0, a_1\in M\), let
\[R_{(a_0, a_1)}=\{(a'_0, a'_1)\in M^{2} :
\exists w\in W
(w(a_0)=a'_0\wedge  w(a_1)=a'_1)\}.\]

We let \[\mm=
\la M,
(R_{(a_0, a_1)})_{(a_0, a_1)\in M^{2}}\ra\]
\[\ma=
\la\mm,a_0\ra\]
\[\mb=
\la \mm,b_0\ra\]
and show that \(\EFG{\goy}\) is non-determined.


The reduction of the language of \(\ma\) and \(\mb\) to one binary
predicate is easy. One just adds a copy of \(\go_3\), together with
its ordering, and a copy of $M\times M$ to the structures with the
projection maps.  Then fix a bijection $\gf$ from $\go_3$ to  $M^2$.
Add a new binary predicate $R$ to the language and interpret $R$ to be
contained in $\go_3 \times M^2$ such that  $R(\gb, (a, b))$ holds if
and only if $R_{\gf(\ga)}(a, b)$ holds. We can now dispense with the
old binary predicates. We have replaced our structure by one in a
finite language without making any difference to who wins the game
\(\EFG{\goy}\).  The extra step of reducing to a single binary
predicate is standard.

An important property of these models is that if \(\ga\in
S\cap S^3_1\),
then \(G_{\ga}\rest M_{\ga}\) is an automorphism of the
restriction
of \(\mm\) to \(M_{\ga}\) and takes \(a_0\) to \(b_0\).
\begin{claim}
\(\I\) does not  win \(\EFG{\goy}\).
\end{claim}

Suppose \(\I\) has a winning strategy \(\gt\).
Again, there is a quick argument which uses CH: Find
\(\ga\in S\) such that \(M_{\ga}\) is closed under
\(\gt\) and \(\cf (\ga)=\goy\). Now \(C_{\ga}\) is cub
on \(\ga\), whence \(G_{\ga}\) maps \(M_{\ga}\) onto
itself. Using \(G_{\ga}\) player \(\II\) can easily
beat \(\gt\), a contradiction.

In the following longer argument we need not assume CH.
Let
\(\gk\) be a large regular cardinal.
 Let \(\mh\) be the expansion of
\(\la H(\gk),\in\ra\) obtained by adding the following
structure to it:
\begin{description}
\item[(H1)] The function \(\ga\mapsto M_{\ga}\).
\item[(H2)] The function \(\ga\mapsto G_{\ga}\).
\item[(H3)] The function \(\ga\mapsto C_{\ga}\).
\item[(H4)] A well-ordering \(<^*\) of the universe.
\item[(H5)] The winning strategy \(\gt\).
\item[(H6)] The sets \(S\) and \(T\).
\end{description}
Let \(\mn=\la N,\in,\ldots \ra\) be an elementary submodel of
\(\mh\)
such that \(\ga=N\cap\go_3\in S\cap S^3_1\).

Now \(C_{\ga}\) is a cub of order-type \(\goy\) on
\(\ga\) and \(G_{\ga}\) maps \(M_{\ga}\)
onto \(M_{\ga}\).
Moreover, \(G_{\ga}\) is a partial isomorphism
from \(\ma\) into \(\mb\).  Provided that \(\gt\)
does not lead \(\I\) to play his moves outside \(M_{\ga}\),
 \(\II\) has on obvious
strategy: he lets \(G_{\ga}\) determine his moves.
So let us assume a sequence
\(\la (x_{\gx},y_{\gx}) :  \gx<\gg\ra\) has been played inside
\(M_{\ga}\) and
\(\gg<\goy\). Let
\(\gb\in C_{\ga}\) such that \(M_{\gb}\) contains the
elements \(x_{\gx},y_{\gx}\) for \(\gx<\gg\). The sequence
\(\la y_{\gx} : \gx<\gg\ra\) is totally determined by
\(G_{\gb}\) and \(\gt\). Since \(G_{\gb}\in N\),
\(\la y_{\gx} : \gx<\gg\ra\in N\), and we are done.

\begin{claim}
\(\II\) does not win \(\EFG{\goy}\).
\end{claim}

Suppose \(\II\) has a winning strategy \(\gt\).
Let \(\mh\) be as above and
\(\mn=\la N,\in,\ldots \ra\) be an elementary
submodel of \(\mh\) such that \(\ga=N\cap\go_3\in T\).
We let \(\I\) play during the first \(\go\) moves of
\(\EFG{\goy}\) a sequence \(\la a_n : n<\go\ra\) in \(\ma\)
such that if
 \(\ga_n\) is the least \(\ga_n\) with
\(a_n\in M_{\ga_n}\), then the sequence
\(\la\ga_n : n<\go\ra\) is ascending and
 \(\sup \{\ga_n : n<\go\}=\ga\).
Let \(\II\) respond following \(\gt\) with
\(\la b_n : n<\go\ra\).
As his move number \(\go\) player \(\I\) plays some element
\(a_{\go}\in M\setminus M_{\ga}\) in \(\ma\) and
\(\II\) answers according to \(\gt\) with \(b_{\go}\).


        For all $i \leq \go$, $R_{(a_0, a_i)}(a_0, a_i)$ holds. Hence
$R_{(a_0, a_i)}(b_0, b_i)$ holds. So there is $w_i$ such that
$w_i(a_0) = b_0$ and $w_i(a_i) = b_i$. Since $G_{\gb_0}(a_0) = b_0$,
by (M5), for each $i$ there is $\gb_i$ so that $G_{\gb_i}(a_i) =
b_i$. We can assume that $\gb_i$ is chosen to be minimal. Notice that
for all $i$, $\gb_i > \ga_i$ and for $i < \go$, $\gb_i \in \mn$. So
$\sup  \set{\gb_i}{i < \go} = \ga$.

        Also, by the same reasoning as above, for each $i < \go$,
$R_{(a_i, a_\go)}(b_i, b_\go)$ holds. Applying (M5), we get that
$G_{\gb_\go}(a_i) = b_i$. Using (M5) again and the minmality of
$\gb_i$, for all $i < \go$, $\gb_i \in C_{\gb_\go}$. Thus \(\ga\) is
a limit of elements of \(C_{\gb_{\go}}\), contradicting \(\ga\in T\).
\fin

\end{document}